\documentclass[oneside,english]{amsart}
\usepackage[T1]{fontenc}
\usepackage[utf8]{inputenc}
\usepackage{prettyref}
\usepackage{amstext}
\usepackage{amsthm}
\usepackage{amssymb}
\usepackage[numbers]{natbib}
\usepackage[all]{xy}

\makeatletter
\numberwithin{equation}{section}
\numberwithin{figure}{section}
\theoremstyle{plain}
\newtheorem{thm}{\protect\theoremname}
\theoremstyle{definition}
\newtheorem{example}[thm]{\protect\examplename}
\theoremstyle{definition}
\newtheorem{defn}[thm]{\protect\definitionname}
\theoremstyle{plain}
\newtheorem{lem}[thm]{\protect\lemmaname}
\theoremstyle{plain}
\newtheorem{cor}[thm]{\protect\corollaryname}
\theoremstyle{remark}
\newtheorem{rem}[thm]{\protect\remarkname}
\theoremstyle{remark}
\newtheorem*{claim*}{\protect\claimname}

\usepackage{mathtools}
\usepackage{braket}

\global\long\def\ns#1{\prescript{\ast}{}{#1}}
\global\long\def\st#1{\prescript{\circ}{}{#1}}

\usepackage[
    type={CC},
    modifier={by-nc-nd},
    version={4.0},
]{doclicense}
\usepackage{hyperxmp}

\makeatother

\usepackage{babel}
\providecommand{\claimname}{Claim}
\providecommand{\corollaryname}{Corollary}
\providecommand{\definitionname}{Definition}
\providecommand{\examplename}{Example}
\providecommand{\lemmaname}{Lemma}
\providecommand{\remarkname}{Remark}
\providecommand{\theoremname}{Theorem}

\begin{document}
\title{A nonstandard construction of direct limit group actions}
\thanks{\copyright~2022 Mathematical Communications. This manuscript version is made available under the \doclicenseLongNameRef.}
\author{Takuma Imamura}
\address{Research Institute for Mathematical Sciences\\
Kyoto University\\
Kitashirakawa Oiwake-cho, Sakyo-ku, Kyoto 606-8502, JAPAN}
\email{timamura@kurims.kyoto-u.ac.jp}
\begin{abstract}
Manevitz and Weinberger (1996) proved that the existence of effective
$K$-Lipschitz $\mathbb{Z}/n\mathbb{Z}$-actions implies the existence
of effective $K$-Lipschitz $\mathbb{Q}/\mathbb{Z}$-actions for all
compact connected manifolds with metrics, where $K$ is a fixed Lipschitz
constant. The $\mathbb{Q}/\mathbb{Z}$-actions were constructed from
suitable actions of a sufficiently large hyperfinite cyclic group
$\ns{\mathbb{Z}}/\gamma\ns{\mathbb{Z}}$ in the sense of nonstandard
analysis. By modifying their construction, we prove that for every
direct system $\left(\Lambda,G_{\lambda},i_{\lambda\mu}\right)$ of
torsion groups with monomorphisms, the existence of effective $K$-Lipschitz
$G_{\lambda}$-actions implies the existence of effective $K$-Lipschitz
$\varinjlim G_{\lambda}$-actions. This generalises Manevitz and Weinberger's
result.
\end{abstract}

\thanks{This work was partly supported by the Morikazu Ishihara (Shikata)
Research Encouragement Fund and by JST ERATO HASUO Metamathematics
for Systems Design Project (No. JPMJER1603).}
\maketitle

\section{\label{sec:Introduction}Introduction}

Let $M$ be a compact connected manifold with a metric. Using nonstandard
analysis, \citet{MW96} proved that if $M$ admits an effective $K$-Lipschitz
action of the cyclic group $\mathbb{Z}/n\mathbb{Z}$ for each $n\in\mathbb{Z}_{+}$,
then $M$ also admits an effective $K$-Lipschitz action of the rational
circle group $\mathbb{Q}/\mathbb{Z}$. (The standard approach to this
type of result with applications can be found in \citet{Wei87}.)
A sketch of the proof is as follows: let $\gamma$ be an infinite
hyperinteger that is divisible by all non-zero integers (e.g. the
factorial $\omega!$ of an arbitrary positive infinite hyperinteger
$\omega$). Then $\mathbb{Q}/\mathbb{Z}$ can be embedded into the
hyperfinite cyclic group $\prescript{\ast}{}{\mathbb{Z}}/\gamma\prescript{\ast}{}{\mathbb{Z}}$
by identifying $k/n+\mathbb{Z}\in\mathbb{Q}/\mathbb{Z}$ with $k\left(\gamma/n\right)+\gamma\ns{\mathbb{Z}}\in\ns{\mathbb{Z}}/\gamma\ns{\mathbb{Z}}$.
By the transfer principle, the nonstandard extension $\ns{M}$ admits
an internal effective $K$-Lipschitz action of $\prescript{\ast}{}{\mathbb{Z}}/\gamma\prescript{\ast}{}{\mathbb{Z}}$.
By restricting the domain and by taking its standard part, we obtain
the desired $\mathbb{Q}/\mathbb{Z}$-action on $M$. The effectiveness
of the resulting action follows from Newman's theorem in the version
of \citet[III.9.6 Corollary]{Bre72}. Their proof requires no advanced
knowledge of transformation group theory. However, their proof contains
an error involving the use of the downward transfer principle (see
the footnote in the proof of \prettyref{thm:Injective-limit-of-torsion-groups}).
Fortunately, their proof can be corrected, as we shall see below.

Considering that $\mathbb{Q}/\mathbb{Z}$ is isomorphic to the direct
limit of $\mathbb{Z}/n\mathbb{Z}\ \left(n\in\mathbb{Z}_{+}\right)$,
it is natural to attempt to generalise Manevitz and Weinberger's result
to direct limits of a more general class of groups. In \prettyref{sec:Preliminaries}
we recall some results from nonstandard analysis and topology. In
\prettyref{sec:Main-results}, we prove that for every direct system
$\mathbf{G}:=\left(\Lambda,G_{\lambda},i_{\lambda\mu}\right)$ of
torsion groups with monomorphisms, if $M$ admits an effective $K$-Lipschitz
$G_{\lambda}$-action for each $\lambda\in\Lambda$, then $M$ also
admits an effective $K$-Lipschitz $\varinjlim\mathbf{G}$-action.
The result on $\mathbb{Q}/\mathbb{Z}$-actions is an immediate corollary
to our result. One can also obtain the following corollary on $\mathbb{Z}\left[1/p\right]/\mathbb{Z}$-actions:
if the cyclic groups $\mathbb{Z}/p^{n}\mathbb{Z}\ \left(n\in\mathbb{N}\right)$
effectively act on $M$ by $K$-Lipschitz maps, then the $p$-Pr\"ufer
group $\mathbb{Z}\left[1/p\right]/\mathbb{Z}$ does as well. In \prettyref{sec:Conclusion},
we conclude the paper by mentioning related works involving nonstandard
approximations of direct and inverse limits.

\section{\label{sec:Preliminaries}Preliminaries}

First of all, we recall the model-theoretic framework of nonstandard
analysis (NSA). We refer to \citet{Rob66,CK90,LW15} for model-theoretic
NSA and \citet{KR04} for axiomatic NSA. The reader is assumed to
be familiar with the rudiments of mathematical logic. NSA uses the
following two universes:
\begin{enumerate}
\item The standard universe $\left(\mathbb{U},\in\right)$. Assume the following:
\begin{description}
\item [{Transitivity}] The underlying set $\mathbb{U}$ is a transitive
set, i.e. $x\in\mathbb{U}$ implies $x\subseteq\mathbb{U}$.
\item [{Richness}] All standard mathematical objects we need (such as groups
and manifolds appear in this paper) belong to $\mathbb{U}$.
\item [{Absoluteness}] All (but finitely many) set-theoretic formulae we
need are absolute with respect to $\mathbb{U}$. In other words, given
a set-theoretic formula $\varphi\left(\vec{x}\right)$ that appears
in this paper (such as ``$U$ is an open set of $X$'' and ``$f$
is continuous at $x$'') and parameters $\vec{a}\in\mathbb{U}$, the
sentence $\varphi\left(\vec{a}\right)$ is true in $\mathbb{U}$ (by
interpreting $\forall$ and $\exists$ as quantifiers over $\mathbb{U}$)
if and only if $\varphi\left(\vec{a}\right)$ is actually true (by
interpreting $\forall$ and $\exists$ as quantifiers over all mathematical
objects). In particular, all axioms of ZFC we need are true in $\mathbb{U}$.
\end{description}
While the existence of such a universe $\mathbb{U}$ is provable in
ZFC by the reflection principle, the reader familiar with category
theory may consider $\mathbb{U}$ as a Grothendieck universe. (The
absoluteness holds for all bounded formulae in this case.) A more
clever approach can be found in \citet{Fef69}.
\item The nonstandard universe $\left(\ns{\mathbb{U}},\ns{\in}\right)$
with the embedding $\ns{\left(-\right)}\colon\mathbb{U}\hookrightarrow\ns{\mathbb{U}}$
satisfying the following principles:
\begin{description}
\item [{Transfer}] For any sentence $\varphi\left(\vec{a}\right)$ with
parameters $\vec{a}$ in $\mathbb{U}$, $\varphi\left(\vec{a}\right)$
is true in $\mathbb{U}$ if and only if $\varphi\left(\ns{\vec{a}}\right)$
is true in $\ns{\mathbb{U}}$. The ``only if'' part is referred to
as \emph{upward transfer}. The ``if'' part is referred to as \emph{downward
transfer}.
\item [{$\kappa$-saturation}] Let $\kappa$ be a \emph{fixed} infinite
cardinal. Let $p\left(\vec{x}\right)$ be a set of formulae with variables
$\vec{x}$ and parameters in $\ns{\mathbb{U}}$. Suppose that $\left|p\left(\vec{x}\right)\right|<\kappa$.
If every finite subset $q\left(\vec{x}\right)$ of $p\left(\vec{x}\right)$
has a solution in $\ns{\mathbb{U}}$, the whole $p\left(\vec{x}\right)$
has a solution in $\ns{\mathbb{U}}$. The limitation of the cardinality
of $p\left(\vec{x}\right)$ cannot be relaxed, because the unlimited
saturation principle leads a contradiction. To prove the main results
of this paper, we only need to assume the following weaker principle.
\item [{Weak saturation}] Let $p\left(\vec{x}\right)$ be a set of formulae
with variables $\vec{x}$ and parameters in $\ns{\mathbb{U}}$. Suppose
that each parameter of $p\left(\vec{x}\right)$ belongs to the image
of the embedding $\ns{\left(-\right)}$. If every finite subset $q\left(\vec{x}\right)$
of $p\left(\vec{x}\right)$ has a solution in $\ns{\mathbb{U}}$,
the whole $p\left(\vec{x}\right)$ has a solution in $\ns{\mathbb{U}}$.
\end{description}
See \citet{CK90} for the construction of $\ns{\mathbb{U}}$.
\end{enumerate}
A mathematical object is said to be \emph{standard} (or \emph{$\mathbb{U}$-small}
in terminology of category theory) if it is an element of $\mathbb{U}$;
\emph{internal} if it is an element of $\ns{\mathbb{U}}$; and \emph{external}
if it is not internal. Given a concept $X$ on $\mathbb{U}$ defined
by a formula $\varphi\left(\vec{x},\vec{a}\right)$, the concept on
$\ns{\mathbb{U}}$ defined by the associated formula $\varphi\left(\vec{x},\ns{\vec{a}}\right)$
is called \emph{internal $X$}, \emph{hyper $X$,} and {*}$X$. We
drop the star $\ns{\left(-\right)}$ unless there is a risk of confusion.
In particular, we identify the {*}membership relation $\ns{\in}$
with the genuine membership relation $\in$.
\begin{example}[The ordered field of hyperreals]
Let $K$ be an ordered field. We may assume without loss of generality
that $K$ is an extension of $\mathbb{Q}$. An element of $K$ is
called an \emph{infinite} (with respect to $\mathbb{Q}$) if its absolute
value is an upper bound of $\mathbb{Q}$. An element of $K$ is called
an \emph{infinitesimal} if its absolute value is a lower bound of
$\mathbb{Q}_{+}$. The ordered field $K$ is said to be \emph{non-Archimedean}
if one of the following equivalent conditions holds: (i) it has an
infinite; (ii) it has a non-zero infinitesimal.

The property that $\mathbb{R}$ is an ordered field can be described
as a formula. By transfer, $\ns{\mathbb{R}}$ is an ordered field.
The field $\ns{\mathbb{R}}$ and its elements are respectively called
the \emph{hyperreal field} and \emph{hyperreal numbers} by the above
convention. The restriction of $\mathbb{U}\hookrightarrow\ns{\mathbb{U}}$
gives an embedding of $\mathbb{R}$ into $\ns{\mathbb{R}}$. Consider
the set $p\left(x\right)=\set{\text{``}x\in\ns{\mathbb{R}}\text{''}}\cup\set{\text{``}a<\left|x\right|\text{''}|a\in\mathbb{Q}}$
of formulae with one variable $x$ (and parameters $\ns{\mathbb{R}}$
and $a\in\mathbb{Q}$). Every finite subset of $p\left(x\right)$
is solvable in $\ns{\mathbb{U}}$, so $p\left(x\right)$ is solvable
in $\ns{\mathbb{U}}$ by weak saturation. The solutions of $p\left(x\right)$
are precisely hyperreal numbers whose absolute values are greater
than all (standard) rational numbers, i.e. infinites. Similarly, one
can obtain non-zero infinitesimals by considering the set $q\left(x\right)=\set{\text{``}x\in\ns{\mathbb{R}}\text{''}}\cup\set{\text{``}x\neq0\text{''}}\cup\set{\text{``}\left|x\right|<a\text{''}|a\in\mathbb{Q}_{+}}$.
Hence $\ns{\mathbb{R}}$ is a non-Archimedean ordered field.
\end{example}

We recall some fundamental results from nonstandard topology.
\begin{defn}[\citet{Rob66}]
Let $\left(X,\tau_{X}\right)$ be a standard topological space. For
$x\in X$, the set $\mu_{X}\left(x\right):=\bigcap_{x\in U\in\tau_{X}}\ns{U}$
is called the \emph{monad} of $x$.
\end{defn}

\begin{defn}[\citet{Rob66}]
Let $X$ be a standard metric space. For $x,y\in\ns{X}$, we say
that $x$ and $y$ are \emph{infinitely close} ($x\approx_{X}y$)
if the {*}distance $\ns{d_{X}}\left(x,y\right)$ is an infinitesimal.
\end{defn}

For a standard metric space $X$, the monad $\mu_{X}\left(x\right)$
of $x\in X$ is precisely the set of all points of $\ns{X}$ infinitely
close to $x$.
\begin{lem}[\citet{Rob66}]
\label{lem:Approximation}Let $X$ be a standard topological space
and $x\in X$. There exists a {*}open set $U$ (i.e. a member of $\ns{\tau_{X}}$)
such that $x\in U\subseteq\mu_{X}\left(x\right)$.
\end{lem}

\begin{proof}
Apply weak saturation to the set 
\[
p\left(U\right):=\set{\text{``}x\in U\text{''},\text{``}U\in\ns{\tau_{X}}\text{''}}\cup\set{\text{``}U\subseteq\ns{V}\text{''}|x\in V\in\tau_{X}}.\qedhere
\]
\end{proof}
\begin{thm}[\citet{Rob66}]
\label{thm:Nbhd}Let $X$ be a standard topological space and $x\in X$.
A subset $U$ of $X$ is a neighbourhood of $x$ if and only if $\mu_{X}\left(x\right)\subseteq\ns{U}$.
\end{thm}

\begin{proof}
The ``only if'' part is trivial by the definition of $\mu_{X}$. To
prove the ``if'' part, suppose that $\mu_{X}\left(x\right)\subseteq\ns{U}$.
By \prettyref{lem:Approximation}, there exists a $V\in\ns{\tau_{X}}$
such that $x\in V\subseteq\mu_{X}\left(x\right)\subseteq\ns{U}$,
i.e. $\ns{U}$ is a {*}neighbourhood of $x$. By downward transfer,
$U$ is a neighbourhood of $x$.
\end{proof}
\begin{cor}[\citet{Rob66}]
\label{cor:Open-sets}Let $X$ be a standard topological space. A
subset $U$ of $X$ is an open set if and only if for all $x\in U$,
$\mu_{X}\left(x\right)\subseteq\ns{U}$.
\end{cor}

\begin{cor}[\citet{Rob66}]
\label{cor:Closed-sets}Let $X$ be a standard topological space.
A subset $F$ of $X$ is a closed set if and only if $\mu_{X}\left(x\right)\cap\ns{F}\neq\varnothing$
implies $x\in F$ for all $x\in X$.
\end{cor}

\begin{thm}[\citet{Rob66}]
A standard topological space $X$ is Hausdorff if and only if $\mu_{X}\left(x\right)\cap\mu_{X}\left(y\right)=\varnothing$
for all distinct $x,y\in X$.
\end{thm}

\begin{proof}
Let $x,y\in X$. It suffices to show that $x$ and $y$ are separable
by neighbourhoods if and only if $\mu_{X}\left(x\right)\cap\mu_{X}\left(y\right)=\varnothing$.
Suppose $x$ and $y$ are separable by neighbourhoods $U_{x}$ and
$U_{y}$. By \prettyref{thm:Nbhd}, $\mu_{X}\left(x\right)\cap\mu_{X}\left(y\right)\subseteq\ns{U_{x}}\cap\ns{U_{y}}=\ns{\left(U_{x}\cap U_{y}\right)}=\ns{\varnothing}=\varnothing$.
Conversely, suppose $\mu_{X}\left(x\right)\cap\mu_{X}\left(y\right)=\varnothing$.
There exist $U_{x},U_{y}\in\ns{\tau_{X}}$ such that $x\in U_{x}\subseteq\mu_{X}\left(x\right)$
and $y\in U_{y}\subseteq\mu_{X}\left(x\right)$ by \prettyref{lem:Approximation}.
Note that $U_{x}\cap U_{y}\subseteq\mu_{X}\left(x\right)\cap\mu_{Y}\left(y\right)=\varnothing$.
By downward transfer, there exist (standard) $U_{x},U_{y}\in\tau_{X}$
such that $x\in U_{x}$, $y\in U_{y}$ and $U_{x}\cap U_{y}=\varnothing$.
\end{proof}
\begin{thm}[\citet{Rob66}]
A standard topological space $X$ is compact if and only if $\ns{X}=\bigcup_{x\in X}\mu_{X}\left(x\right)$.
\end{thm}

\begin{proof}
Suppose $X$ is compact. Let $x\in\ns{X}$. Consider the family of
closed sets of $X$ defined by
\[
\mathcal{F}:=\set{F\subseteq X|F\text{ is closed and }x\in\ns{F}}.
\]
For each $F_{1},\ldots,F_{n}\in\mathcal{F}$, since $x\in\ns{F_{1}}\cap\cdots\cap\ns{F_{n}}\neq\varnothing$,
$F_{1}\cap\cdots\cap F_{n}\neq\varnothing$ by downward transfer.
The intersection $\bigcap\mathcal{F}$ has an element $y$ by the
compactness of $X$. Let $U$ be an arbitrary open neighbourhood of
$y$. Then $X\setminus U\notin\mathcal{F}$, i.e. $x\notin\ns{\left(X\setminus U\right)}=\ns{X}\setminus\ns{U}$.
Hence $x\in\mu_{X}\left(x\right)$, because $U$ was arbitrary.

Suppose $\ns{X}=\bigcup_{x\in X}\mu_{X}\left(x\right)$. Let $\set{F_{i}}_{i\in I}$
be a family of closed subsets with the finite intersection property.
The intersection $\bigcap_{i\in I}\ns{F_{i}}$ has an element $x\in\ns{X}$
by weak saturation. Choose a $y\in X$ such that $x\in\mu_{X}\left(y\right)$.
For each $i\in I$, since $x\in\mu_{X}\left(y\right)\cap\ns{F_{i}}\neq\varnothing$,
$y\in F_{i}$ by \prettyref{cor:Closed-sets}. Therefore $\bigcap_{i\in I}F_{i}$
is non-empty.
\end{proof}
\begin{cor}[\citet{Rob66}]
\label{cor:standard-part-map}If $X$ is a standard compact Hausdorff
space, there exists a unique map $\st{\left(-\right)}\colon\ns{X}\to X$
(called the \emph{standard part map}) such that $x\in\mu_{X}\left(\st{x}\right)$.
\end{cor}

\begin{thm}[\citet{Rob66}]
\label{thm:continuity}A standard map $f\colon X\to Y$ between topological
spaces is continuous at $x\in X$ if and only if $\ns{f}\left[\mu_{X}\left(x\right)\right]\subseteq\mu_{Y}\left(f\left(x\right)\right)$.
\end{thm}

\begin{proof}
Suppose $f$ is continuous at $x$. Let $U$ be a neighbourhood of
$f\left(x\right)$. Then $f^{-1}\left[U\right]$ is a neighbourhood
of $x$. By \prettyref{thm:Nbhd}, $\ns{f}\left[\mu_{X}\left(x\right)\right]\subseteq\ns{f}\left[\ns{f}^{-1}\left[\ns{U}\right]\right]\subseteq\ns{U}$.
Since $U$ was arbitrary, $\ns{f}\left[\mu_{X}\left(x\right)\right]\subseteq\mu_{Y}\left(f\left(x\right)\right)$.

Conversely, suppose $\ns{f}\left[\mu_{X}\left(x\right)\right]\subseteq\mu_{Y}\left(f\left(x\right)\right)$.
Let $U$ be a neighbourhood of $f\left(x\right)$. By \prettyref{lem:Approximation},
there exists a $V\in\ns{\tau_{X}}$ such that $x\in V\subseteq\mu_{X}\left(x\right)$.
Then $\ns{f}\left[V\right]\subseteq\ns{f}\left[\mu_{X}\left(x\right)\right]\subseteq\mu_{Y}\left(f\left(x\right)\right)\subseteq\ns{U}$
by \prettyref{thm:Nbhd}. By downward transfer, there exists a (standard)
$V\in\tau_{X}$ such that $f\left[V\right]\subseteq U$. Hence $f$
is continuous at $x$.
\end{proof}

\section{\label{sec:Main-results}Main results}

\subsection{Nonstandard approximations of direct limits}

Let $\mathbf{G}:=\left(\Lambda,G_{\lambda},i_{\lambda\mu}\right)$
be a standard direct system of groups and homomorphisms. A \emph{cocone}
over $\mathbf{G}$ consists of a group $G$ and a homomorphism $j_{\lambda}\colon G_{\lambda}\to G$
for each $\lambda\in\Lambda$ which makes the following diagram commutative:
\[
\xymatrix{ & G\\
G_{\lambda}\ar[ur]^{j_{\lambda}}\ar[rr]^{i_{\lambda\mu}} &  & G_{\mu}\ar[ul]_{j_{\mu}}
}
\]
Given two cocones $\left(G,j_{\lambda}\right)$ and $\left(H,k_{\lambda}\right)$
over $\mathbf{G}$, a \emph{morphism} between them is a homomorphism
$f\colon G\to H$ (of groups) such that the diagram
\[
\xymatrix{G\ar[rr]^{f} &  & H\\
 & G_{\lambda}\ar[ul]^{j_{\lambda}}\ar[ur]_{k_{\lambda}}
}
\]
is commutative for all $\lambda\in\Lambda$. The collection of cocones
over $\mathbf{G}$ and their morphisms forms a category. The initial
object of this category is called the \emph{colimiting cocone} over
$\mathbf{G}$ and is denoted by $i_{\lambda}\colon G_{\lambda}\to\varinjlim\mathbf{G}$.
The group $\varinjlim\mathbf{G}$ is called the \emph{direct limit}
of $\mathbf{G}$. The direct limit $\varinjlim\mathbf{G}$ can be
constructed as the quotient of the disjoint union $\bigsqcup_{\lambda\in\Lambda}G_{\lambda}:=\bigcup_{\lambda\in\Lambda}\left(\set{\lambda}\times G_{\lambda}\right)$
modulo the equivalence relation $\equiv$ defined by $\left(\lambda,g\right)\equiv\left(\mu,h\right)$
if and only if $i_{\lambda\nu}\left(g\right)=i_{\mu\nu}\left(h\right)$
for some $\nu\geq\lambda,\mu$.
\begin{lem}
\label{lem:infinite-index}There exists an index $\gamma\in\ns{\Lambda}$
such that $\Lambda\leq\gamma$.
\end{lem}

\begin{proof}
Consider the set $p\left(x\right):=\set{\text{``}x\in\Lambda\text{''}}\cup\set{\text{``}\lambda\leq\gamma\text{''}|\lambda\in\Lambda}$
and apply weak saturation.
\end{proof}
\begin{thm}
\label{thm:Embedding-Lemma}Let $\gamma$ be as in \prettyref{lem:infinite-index}.
There exists an embedding $j\colon\varinjlim\mathbf{G}\to\ns{G_{\gamma}}$
such that the following diagram is commutative for every $\lambda\in\Lambda$:
\[
\xymatrix{\varinjlim\mathbf{G}\ar[rr]^{j} &  & \ns{G_{\gamma}}\\
 & G_{\lambda}\ar[ul]^{i_{\lambda}}\ar[ur]_{\ns{i_{\lambda\gamma}}}
}
\]
\end{thm}

\begin{proof}
The homomorphisms $j_{\lambda}:=\ns{i_{\lambda\gamma}}\restriction G_{\lambda}\colon G_{\lambda}\to\ns{G_{\gamma}}\ \left(\lambda\in\Lambda\right)$
form a cocone over $\mathbf{G}$: $j_{\mu}\circ i_{\lambda\mu}=\ns{i_{\mu\gamma}}\restriction G_{\mu}\circ\ns{i_{\lambda\mu}}\restriction G_{\lambda}=\ns{i_{\lambda\gamma}}\restriction G_{\lambda}=j_{\lambda}$.
By the universal mapping property, there exists a (unique) homomorphism
$j\colon\varinjlim\mathbf{G}\to\ns{G_{\gamma}}$ such that the above
diagram is commutative. More specifically, given $g\in\varinjlim\mathbf{G}$,
find a $\lambda\in\Lambda$ and a $g_{\lambda}\in G_{\lambda}$ such
that $g=i_{\lambda}\left(g_{\lambda}\right)$, and then define $j\left(g\right):=\ns{i_{\lambda\gamma}}\left(g_{\lambda}\right)$.

To prove the injectivity, let $g\in\ker j$. Choose $g_{\lambda}\in G_{\lambda}$
such that $g=i_{\lambda}\left(g_{\lambda}\right)$. Then, $\ns{i_{\lambda\gamma}}\left(g_{\lambda}\right)=j\left(g\right)=e$
by definition. Hence there exists a $\mu\in\ns{\Lambda}$ such that
$\ns{i_{\lambda\mu}}\left(g_{\lambda}\right)=e$. By downward transfer,
there exists a $\mu\in\Lambda$ such that $i_{\lambda\mu}\left(g_{\lambda}\right)=e$.
Therefore $g=e$.
\end{proof}
As a corollary, we obtain another construction of direct limits.
\begin{cor}
$\varinjlim\mathbf{G}\cong j\left[\varinjlim\mathbf{G}\right]=\bigcup_{\lambda\in\Lambda}\ns{i_{\lambda\gamma}}\left[G_{\lambda}\right]$.
\end{cor}

The above argument also applies to any other algebraic systems in
the sense of universal algebra which include rings, gyrogroups, lattices
and Heyting algebras.
\begin{example}[\citet{MW96}]
\label{exa:QmodZ}Let $\Lambda:=\mathbb{Z}_{+}$, $G_{\lambda}:=\mathbb{Z}/\lambda\mathbb{Z}$
and $i_{\lambda\mu}\left(k+\lambda\mathbb{Z}\right):=k\left(\mu/\lambda\right)+\mu\mathbb{Z}$,
where the index set $\mathbb{Z}_{+}$ is ordered by the divisibility
relation. Its direct limit is isomorphic to the rational circle group
$\mathbb{Q}/\mathbb{Z}$. The canonical homomorphism $i_{\lambda}\colon\mathbb{Z}/\lambda\mathbb{Z}\to\mathbb{Q}/\mathbb{Z}$
is given by $k+\lambda\mathbb{Z}\mapsto k/\lambda+\mathbb{Z}$. Let
$\gamma$ be a positive hyperinteger divisible by all non-zero (standard)
integers, e.g. the factorial $\omega!$ of an infinite hypernatural
number $\omega\in\ns{\mathbb{N}}\setminus\mathbb{N}$. The direct
limit $\mathbb{Q}/\mathbb{Z}$ is then embedded into $\ns{\mathbb{Z}}/\gamma\ns{\mathbb{Z}}$
by $k/n+\mathbb{Z}\mapsto k\left(\gamma/n\right)+\gamma\ns{\mathbb{Z}}$.
\end{example}

\begin{example}
\label{exa:Prufer}Fix a prime number $p$. Let $\Lambda:=\mathbb{N}$,
$G_{\lambda}:=\mathbb{Z}/p^{\lambda}\mathbb{Z}$ and $i_{\lambda\mu}\left(k+p^{\lambda}\mathbb{Z}\right):=p^{\mu-\lambda}k+p^{\mu}\mathbb{Z}$,
where $\mathbb{N}$ is ordered by the usual linear order $\leq$.
Its direct limit, called the \emph{$p$-Pr\"ufer group}, is isomorphic
to $\mathbb{Z}\left[1/p\right]/\mathbb{Z}$. Let $\gamma$ be an infinite
hypernatural number. The direct limit $\mathbb{Z}\left[1/p\right]/\mathbb{Z}$
can be embedded into $\ns{\mathbb{Z}}/p^{\gamma}\ns{\mathbb{Z}}$
by $\sum_{i=0}^{n}a_{i}p^{-i}+\mathbb{Z}\mapsto\sum_{i=0}^{n}a_{i}p^{\gamma-i}+p^{\gamma}\ns{\mathbb{Z}}$.
\end{example}

\begin{example}
Let $k$ be a commutative field. Let $\Lambda:=\mathbb{N}$ (with
$\leq$), $G_{\lambda}:=\mathit{GL}\left(\lambda,k\right)$ and
\[
i_{\lambda\mu}\left(A_{\lambda}\right):=\left(\begin{array}{cc}
A_{\lambda} & O_{\lambda,\mu-\lambda}\\
O_{\mu-\lambda,\lambda} & I_{\mu-\lambda}
\end{array}\right).
\]
The direct limit $\mathit{GL}\left(\infty,k\right)$ is the group
of regular $\infty\times\infty$-matrices of the form:
\[
\left(\begin{array}{cc}
A_{\lambda} & O_{\lambda,\infty}\\
O_{\infty,\lambda} & I_{\infty}
\end{array}\right),\ A_{\lambda}\in\mathit{GL}\left(\lambda,k\right).
\]
Let $\gamma\in\ns{\mathbb{N}}$ be an infinite hypernatural number.
Then $\mathit{GL}\left(\infty,k\right)$ is embedded into $\ns{\mathit{GL}\left(\gamma,\ns{k}\right)}$
by
\[
j\left(\begin{array}{cc}
A_{\lambda} & O_{\lambda,\infty}\\
O_{\infty,\lambda} & I_{\infty}
\end{array}\right):=\left(\begin{array}{cc}
A_{\lambda} & O_{\lambda,\gamma-\lambda}\\
O_{\gamma-\lambda,\lambda} & I_{\gamma-\lambda}
\end{array}\right).
\]
\end{example}

\subsection{Construction of direct limit group actions}
\begin{defn}
Suppose that a group $G$ acts on a metric space $M$. Let $K>0$.
The action is said to be \emph{effective} if for each $g\in G\setminus\set{e_{G}}$,
$d_{M}\left(x,gx\right)>0$ for some $x\in M$. The action is said
to be \emph{$K$-Lipschitz} if $d_{M}\left(gx,gy\right)\leq Kd_{M}\left(x,y\right)$
for all $g\in G$ and all $x,y\in M$.
\end{defn}

Our main theorem is the following.
\begin{thm}
\label{thm:Injective-limit-of-torsion-groups}Let $M$ be a compact
connected manifold with a metric $d_{M}$. Let $\mathbf{G}:=\left(\Lambda,G_{\lambda},i_{\lambda\mu}\right)$
be a direct system of torsion groups, where $i_{\lambda\mu}$ is a
monomorphism for all $\lambda\leq\mu$. If there exists an effective
$K$-Lipschitz $G_{\lambda}$-action on $M$ for each $\lambda\in\Lambda$,
then there exists an effective $K$-Lipschitz $\varinjlim\mathbf{G}$-action
on $M$.
\end{thm}

\begin{rem}
If there exists for each $\lambda\in\Lambda$ an effective $K$-Lipschitz
$G_{\lambda}$-action $\Phi_{\lambda}$ on $M$ such that $\Phi_{\mu}\circ i_{\lambda\mu}=\Phi_{\lambda}$
for all $\lambda\leq\mu$, then there exists an effective $K$-Lipschitz
$\varinjlim\mathbf{G}$-action on $M$. In order to prove it, we just
construct the action by gluing the given actions: $\Psi\left(i_{\lambda}\left(g\right),x\right):=\Phi_{\lambda}\left(g,x\right)$.
We do not assume such a coherency condition in \prettyref{thm:Injective-limit-of-torsion-groups}.
\end{rem}

Before proving the theorem, we consider some direct consequences (see
\prettyref{exa:QmodZ} and \prettyref{exa:Prufer}).
\begin{cor}[\citet{MW96}]
Let $M$ be a compact connected manifold with a metric. If there
exists an effective $K$-Lipschitz $\mathbb{Z}/n\mathbb{Z}$-action
on $M$ for each $n\in\mathbb{Z}_{+}$, then there exists an effective
$K$-Lipschitz $\mathbb{Q}/\mathbb{Z}$-action on $M$.
\end{cor}

\begin{cor}
Let $M$ be a compact connected manifold with a metric and $p$ a
prime number. If there exists an effective $K$-Lipschitz $\mathbb{Z}/p^{n}\mathbb{Z}$-action
on $M$ for each $n\in\mathbb{N}$, then there exists an effective
$K$-Lipschitz $\mathbb{Z}\left[1/p\right]/\mathbb{Z}$-action on
$M$.
\end{cor}

Note that $\mathbb{Q}/\mathbb{Z}$ and $\mathbb{Z}\left[1/p\right]/\mathbb{Z}$
are locally finite, i.e. every finitely generated subgroup is finite.
In fact, the above corollaries are a consequence of a more general
corollary on locally finite group actions.
\begin{cor}
Let $M$ be a compact connected manifold with a metric and let $G$
be a locally finite group. If every finite subgroup $H$ of $G$ acts
effectively on $M$ by $K$-Lipschitz maps, then $G$ acts effectively
on $M$ by $K$-Lipschitz maps.
\end{cor}

\begin{proof}
Consider the set $\Lambda$ of all finite subgroups of $G$ ordered
by the inclusion relation $\subseteq$. We first verify that $\Lambda$
is a directed set. Let $H_{1},\ldots,H_{n}\in\Lambda$. Since $G$
is locally finite, the group $H'$ generated by $H_{1}\cup\cdots\cup H_{n}$
is finite, i.e. $H'\in\Lambda$. The finite group $H'$ is an upper
bound of $\set{H_{1},\ldots,H_{n}}$. For $H,H'\in\Lambda$ with $H\subseteq H'$,
let $i_{HH'}\colon H\to H'$ be the inclusion map. Then $\mathbf{G}:=\left(\Lambda,H,i_{HH'}\right)$
forms a direct system of torsion groups with monomorphisms.

By the local finiteness, the group $\braket{g}$ generated by $g$
is finite for each $g\in G$, so $G=\bigcup_{g\in G}\braket{g}\subseteq\bigcup_{H\in\Lambda}H\subseteq G$.
It is easy to see that the direct limit $\varinjlim\mathbf{G}$ is
isomorphic to $\bigcup_{H\in\Lambda}H$, which is precisely $G$.
The statement of the corollary now follows by \prettyref{thm:Injective-limit-of-torsion-groups}.
\end{proof}
In \prettyref{thm:Injective-limit-of-torsion-groups}, to prove effectiveness,
we employ the following version of Newman's theorem.
\begin{thm}[{\citet[Theorem 2]{Dre69}}]
Let $M$ be a connected manifold with a metric $d_{M}$. There exists
a constant $\varepsilon:=\varepsilon\left(M,d_{M}\right)>0$ such
that for every effective action of a finite group $G$ on $M$, there
exist $g\in G$ and $x\in M$ such that $d_{M}\left(x,gx\right)\geq\varepsilon$.
\end{thm}

A group action $G\curvearrowright M$ is said to be \emph{$\varepsilon$–effective}
if for every $g\in G\setminus\set{e_{G}}$ there exists an $x\in M$
such that the orbit $Gx$ has diametre at least $\varepsilon$. In
this terminology, Newman's theorem states every effective action of
a finite group is $\varepsilon$–effective, where $\varepsilon>0$
depends only on $M$.
\begin{proof}[Proof of \prettyref{thm:Injective-limit-of-torsion-groups}]
By the absoluteness of $\mathbb{U}$, we may assume without loss
of generality that all the objects appeared in the statement (such
as $M$ and $\mathbf{G}$) are standard. For simplicity, denote $G=\varinjlim\mathbf{G}$.
Let $j\colon G\to\ns{G_{\gamma}}$ be the embedding of \prettyref{thm:Embedding-Lemma},
where $\gamma\in\ns{\Lambda}$ is an upper bound of $\Lambda$. By
upward transfer, there exists an internal effective $K$-Lipschitz
action $\Phi\colon\ns{G_{\gamma}}\times\ns{M}\to\ns{M}$. Since $M$
is compact Hausdorff, each point $x\in\prescript{\ast}{}{M}$ is infinitely
close to a unique point $\prescript{\circ}{}{x}\in M$ (see \prettyref{cor:standard-part-map}).
Now define a map $\Psi\colon G\times M\to M$ by putting
\[
\Psi\left(g,x\right):=\st{\left(\Phi\left(j\left(g\right),x\right)\right)}.\phantom{\qedhere}
\]
\begin{claim*}
$\Psi$ is an action.
\end{claim*}
\begin{proof}
Let $g,h\in G$ and $x\in M$. Since $j\colon G\to\ns{G_{\gamma}}$
and $\Phi\colon\ns{G_{\gamma}}\to\ns{\left(\mathrm{Aut}\left(M\right)\right)}$
are homomorphisms, we have that
\begin{align*}
\Psi\left(e,x\right) & =\st{\left(\Phi\left(j\left(e\right),x\right)\right)}\\
 & =\st{\left(\Phi\left(e,x\right)\right)}\\
 & =\st{x}\\
 & =x,
\end{align*}
and

\begin{align*}
\Psi\left(gh,x\right) & =\st{\left(\Phi\left(j\left(gh\right),x\right)\right)}\\
 & =\st{\left(\Phi\left(j\left(g\right)j\left(h\right),x\right)\right)}\\
 & =\st{\left(\Phi\left(j\left(g\right),\Phi\left(j\left(h\right),x\right)\right)\right)}\\
 & =\st{\left(\Phi\left(j\left(g\right),\st{\left(\Phi\left(j\left(h\right),x\right)\right)}\right)\right)}\\
 & =\Psi\left(g,\Psi\left(h,x\right)\right).
\end{align*}
Note that the fourth equality of the latter comes from the $K$-Lipschitz
property of $\Phi$:
\begin{align*}
\MoveEqLeft\ns{d_{M}}\left(\Phi\left(j\left(g\right),\Phi\left(j\left(h\right),x\right)\right),\Phi\left(j\left(g\right),\st{\left(\Phi\left(j\left(h\right),x\right)\right)}\right)\right)\\
 & \leq K\ns{d_{M}}\left(\Phi\left(j\left(h\right),x\right),\st{\left(\Phi\left(j\left(h\right),x\right)\right)}\right)\\
 & =\text{finite}\times\text{infinitesimal}\\
 & =\text{infinitesimal}.
\end{align*}
Hence $\Phi\left(j\left(g\right),\Phi\left(j\left(h\right),x\right)\right)$
and $\Phi\left(j\left(g\right),\st{\left(\Phi\left(j\left(h\right),x\right)\right)}\right)$
have the same standard part.
\end{proof}
\begin{claim*}
$\Psi$ is $K$-Lipschitz.
\end{claim*}
\begin{proof}
Let $g\in G$ and $x,y\in M$. Since $\Psi\left(g,x\right)\approx_{M}\Phi\left(j\left(g\right),x\right)$
and $\Psi\left(g,y\right)\approx_{M}\Phi\left(j\left(g\right),y\right)$,
\[
\left(\Psi\left(g,x\right),\Psi\left(g,y\right)\right)\approx_{M\times M}\left(\Phi\left(j\left(g\right),x\right),\Phi\left(j\left(g\right),y\right)\right).
\]
By the the continuity of the metric function $d_{M}\colon M\times M\to\mathbb{R}$
(\prettyref{thm:continuity}), we have

\begin{align*}
d_{M}\left(\Psi\left(g,x\right),\Psi\left(g,y\right)\right) & \approx_{\mathbb{R}}\ns{d_{M}}\left(\Phi\left(j\left(g\right),x\right),\Phi\left(j\left(g\right),y\right)\right)\\
 & \leq Kd_{M}\left(x,y\right).
\end{align*}
It follows that $d_{M}\left(\Psi\left(g,x\right),\Psi\left(g,y\right)\right)\leq Kd_{M}\left(x,y\right)$.
\end{proof}
\begin{claim*}
$\Psi$ is effective.
\end{claim*}
\begin{proof}
Let $g\in G\setminus\set{e}$. Choose a $g_{\lambda}\in G_{\lambda}$
such that $g=i_{\lambda}\left(g_{\lambda}\right)$. Since $i_{\lambda}$
is a homomorphism, $g_{\lambda}$ is not the unit element. Since $G_{\lambda}$
is a torsion group, the group $\braket{g_{\lambda}}$ generated by
$g_{\lambda}$ is a finite subgroup of $\ns{G_{\lambda}}$. (Note
that $\ns{A}=A$ holds for all standard finite sets $A$ by upward
transfer.) Consider the internal action $\Phi_{\lambda}\colon\braket{g_{\lambda}}\times\ns{M}\to\ns{M}$
defined by
\[
\Phi_{\lambda}\left(h,x\right):=\Phi\left(\ns{i_{\lambda\gamma}}\left(h\right),x\right).
\]
Since $\ns{i_{\lambda\gamma}}$ is injective by upward transfer, $\Phi_{\lambda}$
is effective. By Newman's theorem and upward transfer, there exists
a \emph{standard} constant $\varepsilon:=\varepsilon\left(M,d_{M}\right)>0$
such that ``there exist an $h\in\braket{g_{\lambda}}$ and an $x\in\ns{M}$
such that $\ns{d_{M}}\left(x,\Phi_{\lambda}\left(h,x\right)\right)\geq\varepsilon$''.\footnote{The downward transfer principle cannot be applied to the quoted statement,
because it contains a nonstandard object, namely $\Phi_{\lambda}$.
\citet[p. 152, ll. 21–24]{MW96} accidentally applied the downward
transfer principle to the corresponding statement in the original
proof. As you can see, this error can be avoided.} There exists a \emph{standard} $n\in\mathbb{N}$ such that $h=g_{\lambda}^{n}$.
Then $\ns{d_{M}}\left(x,\Phi_{\lambda}\left(g_{\lambda}^{n},x\right)\right)\geq\varepsilon$
holds. Since $\Phi$ is $K$-Lipschitz,
\begin{align*}
\Psi\left(g^{n},\st{x}\right) & =\st{\left(\Phi\left(j\left(g^{n}\right),\st{x}\right)\right)}\\
 & \approx_{M}\Phi\left(j\left(g^{n}\right),\st{x}\right)\\
 & \approx_{M}\Phi\left(j\left(g^{n}\right),x\right)\\
 & =\Phi\left(\ns{i_{\lambda\gamma}}\left(g_{\lambda}^{n}\right),x\right)\\
 & =\Phi_{\lambda}\left(g_{\lambda}^{n},x\right).
\end{align*}
By the continuity of $d_{M}$,
\[
d_{M}\left(\st{x},\Psi\left(g^{n},\st{x}\right)\right)\approx_{\mathbb{R}}\ns{d_{M}}\left(x,\Phi_{\lambda}\left(g_{\lambda}^{n},x\right)\right)\geq\varepsilon.
\]
Hence $d_{M}\left(\st{x},\Psi\left(g^{n},\st{x}\right)\right)\geq\varepsilon>0$.
Since $\Psi\left(g\right)^{n}=\Psi\left(g^{n}\right)\neq\mathrm{id}_{M}$,
it follows that $\Psi\left(g\right)\neq\mathrm{id}_{M}$.
\end{proof}
\end{proof}

\section{\label{sec:Conclusion}Conclusion}

Our proof can be summarised as follows. For each index $\lambda\in\Lambda$,
there exists an effective $K$-Lipschitz action $G_{\lambda}\curvearrowright M$.
Fix an infinitely large index $\gamma\in\ns{\Lambda}$. By transfer,
there exists an effective $K$-Lipschitz action $\Phi\colon\ns{G_{\gamma}}\curvearrowright\ns{M}$
in $\ns{\mathbb{U}}$. Since the direct limit $\varinjlim\mathbf{G}$
can be embedded into $\ns{G_{\gamma}}$, the desired action is obtained
as the restriction of the standard part $\st{\Phi}\colon\ns{G_{\gamma}}\curvearrowright M$.
The effectiveness of the resulting action follows from Newman's theorem.

The crux of this paper is the idea of approximating categorical limits
by nonstandard objects rather than the results themselves. This enables
us to study categorical limits with nonstandard analysis. Here are
some examples of nonstandard approximations of direct and inverse
limits.

\subsection{\v{C}ech theory and McCord theory}

First, recall the definition of \v{C}ech (co)homology groups following
\citet{Dow52}. Let $X$ be a topological space and $G$ an abelian
group. The family $\mathrm{Cov}_{X}$ of all open covers of $X$ forms
a (downward) directed set with respect to the refinement relation.
If $\lambda$ is a refinement of $\mu$, there exists a (canonical)
homomorphism $V\left(\lambda\right)\to V\left(\mu\right)$ of Vietoris
complexes. Here the \emph{Vietoris complex} $V\left(\lambda\right)$
is the simplicial set where $a_{0},\ldots a_{p}\in X$ span a $p$–simplex
if $a_{0},\ldots a_{p}\in U$ for some $U\in\lambda$. The \emph{\v{C}ech
(co) homology groups} of $X$ with coefficients in $G$ are then defined
as the limits:
\begin{align*}
\check{H}_{\bullet}\left(X;G\right) & :=\varprojlim_{\lambda\in\mathrm{Cov}_{X}}H_{\bullet}\left(V\left(\lambda\right);G\right).\\
\check{H}^{\bullet}\left(X;G\right) & :=\varinjlim_{\lambda\in\mathrm{Cov}_{X}}H^{\bullet}\left(V\left(\lambda\right);G\right).
\end{align*}
By \prettyref{lem:Approximation}, $\check{H}^{\bullet}\left(X\right)$
can be embedded into $\ns{H^{\bullet}}\left(\ns{V}\left(\lambda\right)\right)$
for all infinitely fine $\lambda\in\ns{\mathrm{Cov}_{X}}$. This gives
a nonstandard construction of \v{C}ech cohomology.

\citet{McC72} has given a much deeper construction of \v{C}ech (co)homology.
Let $X$ be a standard topological space and $G$ an internal abelian
group. For $p\in\mathbb{N}$, a $p+1$–tuple $\left(a_{0},\ldots,a_{p}\right)$
from $\ns{X}$ are called a \emph{$p$–microsimplex} if $a_{0},\ldots,a_{p}\in\mu_{X}\left(x\right)$
for some $x\in X$. Denote the set of $p$-microsimplexes by $\Delta^{p}$.
A hyperfinite formal sum $\sum_{i=1}^{n}g_{i}\sigma_{i}$ of $p$–microsimplexes
$\sigma_{i}$ with coefficients $g_{i}\in G$, where $\set{g_{i}}_{i=1}^{n}$
and $\set{\sigma_{i}}_{i=1}^{n}$ are both internal, is called a \emph{$p$–microchain}.
(Formally, a $p$-microchain is an internal map $\sigma\colon\ns{X}^{p+1}\to G$
whose support is a hyperfinite subset of $\Delta^{p}$.) The set $M_{p}\left(X;G\right)$
of $p$-microchains forms an abelian group with respect to the usual
addition. The boundary homomorphisms $\partial_{p}\colon M_{p}\left(X;G\right)\to M_{p-1}\left(X;G\right)$
are defined by
\[
\partial_{p}\sum_{i=1}^{n}\left(a_{0}^{i},\ldots,a_{p}^{i}\right):=\sum_{i=1}^{n}\sum_{j=0}^{p}\left(-1\right)^{j}g_{i}\left(a_{0}^{i},\ldots,a_{j-1}^{i},a_{j+1}^{i},\ldots,a_{p}^{i}\right).
\]
Thus $M_{\bullet}\left(X;G\right)$ forms a chain complex.
\[
\xymatrix{0 & M_{0}\left(X;G\right)\ar[l] & M_{1}\left(X;G\right)\ar[l]_{\partial_{1}} & M_{2}\left(X;G\right)\ar[l]_{\partial_{2}} & \cdots\cdots\ar[l]}
\]
The \emph{McCord homology groups} of $X$ with coefficients in $G$
are defined by
\[
H_{\bullet}^{M}\left(X;G\right):=H_{\bullet}\left(M_{\bullet}\left(X;G\right)\right),
\]
where $H_{\bullet}$ is the homology functor of chain complexes. Roughly
speaking, McCord's homology of $X$ is the homology of the Vietoris
complex of the monads $\set{\mu_{X}\left(x\right)|x\in X}$. Indeed
the following isomorphism results are known.
\begin{thm}[\citet{Gar78}]
Assume that $\ns{\mathbb{U}}$ is sufficiently saturated. Let $X$
be a standard compact space and $G$ a standard abelian group. Then
$H_{\bullet}^{M}\left(X;\ns{G}\right)\cong\check{H}_{\bullet}\left(X;\ns{G}\right)$.
\end{thm}

\begin{thm}[\citet{Kor10b}]
Assume that $\ns{\mathbb{U}}$ is sufficiently saturated. Let $X$
be a standard completely regular space and $G$ a standard abelian
group. Then $H_{\bullet}^{M}\left(X;\ns{G}\right)\cong\varinjlim_{K}\check{H}_{\bullet}\left(K;\ns{G}\right)$,
where the direct limit runs over all compact subspaces of $X$.
\end{thm}

As a result of taking inverse limits, \v{C}ech homology may violate
the exactness axiom in the Eilenberg–Steenrod axioms depending on
the choice of the coefficient group (see \citet{MS82}). \citet{Gar78}
proved that \v{C}ech homology is exact for all compact pairs if and
only if the coefficient group is equationally compact. In contrast,
McCord homology satisfies the exactness axiom for all coefficient
groups (\citet{McC72}). See also \citet{Kor10a}.

One can also consider a cohomological counterpart of McCord's theory.
In contrast with McCord homology, there are at least two different
definitions of McCord cohomology. One is the homology based on \emph{external}
cochains. Let $X$ be a standard topological space and $G$ an abelian
group. Define the cochain complex $M^{\bullet}\left(X;G\right)$ by
\[
M^{p}\left(X;G\right):=\hom\left(M_{p}\left(X;\ns{\mathbb{Z}}\right),G\right),
\]
where the coboundary homomorphisms $d_{p}\colon M^{p}\left(X;G\right)\to M^{p+1}\left(X;G\right)$
are defined as usual:
\[
d_{p}\varphi\left(u\right):=\varphi\left(\partial_{p}u\right).
\]
Note that $\varphi\colon M_{p}\left(X;\ns{\mathbb{Z}}\right)\to G$
may not be determined by its values on $\Delta^{p}$. The \emph{McCord
cohomology groups} of $X$ with coefficients in $G$ are defined as
\[
H_{M}^{\bullet}\left(X;G\right):=H^{\bullet}\left(M^{\bullet}\left(X;G\right)\right),
\]
where $H^{\bullet}$ is the cohomology functor of cochain complexes.
\begin{thm}[\citet{Ziv87}]
Assume that $\ns{\mathbb{U}}$ is sufficiently saturated. Let $X$
be a standard locally contractible paracompact space and $G$ an abelian
group. Then $H_{M}^{\bullet}\left(X;\ns{G}\right)\cong\check{H}^{\bullet}\left(X;\hom\left(\ns{\mathbb{Z}},G\right)\right)$.
\end{thm}

Assume that $G$ is internal. We say that a $p$-cochain $\varphi\in M^{p}\left(X;G\right)$
is \emph{essentially internal} if there exists an internal homomorphism
$f\colon\ns{\left(\mathbb{Z}\braket{X^{p+1}}\right)}\to G$ such that
$\varphi\restriction M_{p}\left(X;\ns{\mathbb{Z}}\right)=f\restriction M_{p}\left(X;\ns{\mathbb{Z}}\right)$,
where $\mathbb{Z}\braket{X^{p+1}}$ denotes the free $\mathbb{Z}$-module
generated by $X^{p+1}$. As a consequence of upward transfer, each
essentially internal cochain $\varphi\colon M_{p}\left(X;\ns{\mathbb{Z}}\right)\to G$
is completely determined by its values on $\Delta^{p}$, so $\varphi$
can be identifies with a map $\varphi\restriction\Delta^{p}\colon\Delta^{p}\to G$
having an extension $f\restriction\ns{X^{p+1}}\colon\ns{X^{p+1}}\to G$
in $\ns{\mathbb{U}}$. The set $M_{\mu}^{\bullet}\left(X;G\right)$
of essentially internal cochains forms a subcomplex of $M^{\bullet}\left(X;G\right)$.
Finally define
\[
H_{\mu}^{\bullet}\left(X;G\right):=H^{\bullet}\left(M_{\mu}^{\bullet}\left(X;G\right)\right).
\]

\begin{thm}[\citet{Ziv87}]
Assume that $\ns{\mathbb{U}}$ is sufficiently saturated. Let $X$
be a standard paracompact space and $G$ an internal abelian group.
Then $H_{\mu}^{\bullet}\left(X;G\right)\cong\check{H}^{\bullet}\left(X;G\right)$.
\end{thm}

The uniform versions of \v{C}ech and McCord theories are studied
in \citet{Ima16,Ima21}.

\subsection{Shape theory}

\citet{Wat78} introduced and studied the envelope functor of metric
spaces, which is a nonstandard analogue of Borsuk's shape theory.
Intuitively, the envelope of a metric space $X$ is the strong homotopy
type of the infinitesimal boldification of $\ns{X}$ within an ambient
normed linear space $\ns{Y}$. Shape theory can be formulated in terms
of inverse systems (see \citet{MS82}). Wattenberg's theory is then
considered as an example of nonstandard approximations of inverse
limits.

\subsection{Ends}

Let $\left(X,\xi\right)$ be a pointed metric space. For each $r>0$,
let $E_{r}$ be the set of all unbounded connected components of $X\setminus B_{r}\left(\xi\right)$,
where $B_{r}\left(\xi\right)$ denotes the open ball. If $r<s$, there
exists a canonical surjection $E_{s}\to E_{r}$ which sends $Q\in E_{s}$
to $Q'\in E_{r}$ so that $Q\subseteq Q'$. The elements of the inverse
limit
\[
e\left(X,\xi\right):=\varprojlim E_{r}
\]
in the category of sets are called \emph{ends} (based at $\xi$).
This notion plays a central role in geometric group theory (see e.g.
\citet{BH99,Loh17}). The nonstandard construction of ends can be
found in \citet{Gol11,Ima20}.

\subsection*{Acknowledgement}

The author is grateful to the anonymous reviewer(s) for their invaluable
comments which improved the manuscript.

\bibliographystyle{IEEEtranSN}
\bibliography{bibliography}

\begin{thebibliography}{23}
\providecommand{\natexlab}[1]{#1}
\providecommand{\url}[1]{#1}
\csname url@samestyle\endcsname
\providecommand{\newblock}{\relax}
\providecommand{\bibinfo}[2]{#2}
\providecommand{\BIBentrySTDinterwordspacing}{\spaceskip=0pt\relax}
\providecommand{\BIBentryALTinterwordstretchfactor}{4}
\providecommand{\BIBentryALTinterwordspacing}{\spaceskip=\fontdimen2\font plus
\BIBentryALTinterwordstretchfactor\fontdimen3\font minus
  \fontdimen4\font\relax}
\providecommand{\BIBforeignlanguage}[2]{{%
\expandafter\ifx\csname l@#1\endcsname\relax
\typeout{** WARNING: IEEEtranSN.bst: No hyphenation pattern has been}%
\typeout{** loaded for the language `#1'. Using the pattern for}%
\typeout{** the default language instead.}%
\else
\language=\csname l@#1\endcsname
\fi
#2}}
\providecommand{\BIBdecl}{\relax}
\BIBdecl

\bibitem[Bredon(1972)]{Bre72}
G.~E. Bredon, \emph{Introduction to {C}ompact {T}ransformation {G}roups}.\hskip
  1em plus 0.5em minus 0.4em\relax Academic Press, 1972.

\bibitem[Bridson and H\"afliger(1999)]{BH99}
M.~R. Bridson and A.~H\"afliger, \emph{Metric {S}paces of {N}on-{P}ositive
  {C}urvature}.\hskip 1em plus 0.5em minus 0.4em\relax Springer, 1999.

\bibitem[Chang and Keisler(1990)]{CK90}
C.~C. Chang and H.~J. Keisler, \emph{Model {T}heory}, 3rd~ed.\hskip 1em plus
  0.5em minus 0.4em\relax North-Holland, 1990.

\bibitem[Dowker(1952)]{Dow52}
C.~H. Dowker, ``Homology {G}roups of {R}elations,'' \emph{Ann. of Math.},
  vol.~56, no.~1, pp. 84--95, 1952.

\bibitem[Dress(1969)]{Dre69}
A.~Dress, ``Newman's theorems on transformation groups,'' \emph{Topology},
  vol.~8, no.~2, pp. 203--207, 1969.

\bibitem[Feferman(1969)]{Fef69}
S.~Feferman, ``Set-theoretic foundations of category theory,'' in \emph{Reports
  of the Midwest Category Seminar III}, S.~MacLane, Ed.\hskip 1em plus 0.5em
  minus 0.4em\relax Springer, 1969, pp. 201--247.

\bibitem[Garavaglia(1978)]{Gar78}
S.~Garavaglia, ``Homology with equationally compact coefficients,'' \emph{Fund.
  Math.}, vol. 100, no.~1, pp. 89--95, 1978.

\bibitem[Goldbring(2011)]{Gol11}
I.~Goldbring, ``Ends of groups: a nonstandard perspective,'' \emph{J. Log.
  Anal.}, vol.~3, no.~7, pp. 1--28, 2011.

\bibitem[Imamura(2016)]{Ima16}
T.~Imamura, ``Nonstandard homology theory for uniform spaces,'' \emph{Topology
  Appl.}, vol. 209, pp. 22--29, 2016, corrigendum in
  DOI:10.13140/RG.2.2.36585.75368.

\bibitem[Imamura(2020)]{Ima20}
------, ``A nonstandard invariant of coarse spaces,'' \emph{Grad. J. Math.},
  vol.~5, no.~1, pp. 1--8, 2020.

\bibitem[Imamura(2021)]{Ima21}
------, ``Relationship among various {V}ietoris-type and microsimplicial
  homology theories,'' \emph{Arch. Math. (Brno)}, vol.~57, no.~3, pp. 131--150,
  2021.

\bibitem[Kanovei and Reeken(2004)]{KR04}
V.~Kanovei and M.~Reeken, \emph{Nonstandard {A}nalysis, {A}xiomatically}.\hskip
  1em plus 0.5em minus 0.4em\relax Springer, 2004.

\bibitem[Korppi(2010{\natexlab{b}})]{Kor10a}
T.~Korppi, ``Vanishing of derived limits of non-standard inverse systems,''
  \emph{Topology Appl.}, vol. 157, no.~17, pp. 2692--2703, 2010.

\bibitem[Korppi(2010{\natexlab{a}})]{Kor10b}
------, ``On the homology of compact spaces by using non-standard methods,''
  \emph{Topology Appl.}, vol. 157, no.~17, pp. 2704--2714, 2010.

\bibitem[Loeb and Wolff(2015)]{LW15}
P.~A. Loeb and M.~P.~H. Wolff, Eds., \emph{Nonstandard {A}nalysis for the
  {W}orking {M}athematician}, 2nd~ed.\hskip 1em plus 0.5em minus 0.4em\relax
  Springer, 2015.

\bibitem[L\"oh(2017)]{Loh17}
C.~L\"oh, \emph{Geometric {G}roup {T}heory}.\hskip 1em plus 0.5em minus
  0.4em\relax Springer, 2017.

\bibitem[Manevitz and Weinberger(1996)]{MW96}
L.~M. Manevitz and S.~Weinberger, ``Discrete circle actions: {A} note using
  non-standard analysis,'' \emph{Israel J. Math.}, vol.~94, no.~1, pp.
  147--155, 1996.

\bibitem[Marde\v{s}i\'{c} and Segal(1982)]{MS82}
S.~Marde\v{s}i\'{c} and J.~Segal, \emph{Shape {T}heory}, ser. North-Holland
  Mathematical Library.\hskip 1em plus 0.5em minus 0.4em\relax Amsterdam-New
  York-Oxford: North-Holland, 1982, vol.~26.

\bibitem[McCord(1972)]{McC72}
M.~C. McCord, ``Non-standard analysis and homology,'' \emph{Fund. Math.},
  vol.~74, no.~1, pp. 21--28, 1972.

\bibitem[Robinson(1966)]{Rob66}
A.~Robinson, \emph{Non-standard {A}nalysis}.\hskip 1em plus 0.5em minus
  0.4em\relax North-Holland, 1966.

\bibitem[Wattenberg(1978)]{Wat78}
F.~Wattenberg, ``Nonstandard analysis and the theory of shape,'' \emph{Fund.
  Math.}, vol.~98, no.~1, pp. 41--60, 1978.

\bibitem[Weinberger(1987)]{Wei87}
S.~Weinberger, ``Free $\mathbb{Q}/\mathbb{Z}$ {A}ctions,'' \emph{Comment. Math.
  Helv.}, vol.~62, pp. 450--464, 1987.

\bibitem[{\v{Z}}ivaljevi\'c(1987)]{Ziv87}
R.~T. {\v{Z}}ivaljevi\'c, ``On a cohomology theory based on hyperfinite sums of
  microsimplexes,'' \emph{Pacific J. Math.}, vol. 128, no.~1, pp. 201--208,
  1987.

\end{thebibliography}

\end{document}